\documentclass{article}
\usepackage{mypaper}

\usepackage[english]{babel}
\usepackage[utf8]{inputenc}
\usepackage[T1]{fontenc} % This can slightly change font appearance

\usepackage{graphicx}
\usepackage{epstopdf}

\usepackage[mode=buildnew]{standalone}
\usepackage{tikz}
\usetikzlibrary{shapes,arrows,positioning,external,calc,math}
\usepackage{pgfplots}
\pgfplotsset{compat=1.16}
\usepackage{url}
\usepackage{verbatim}
\usepackage{multirow, booktabs, adjustbox} % Table imports
\usepackage{algorithm}
\usepackage{algorithmic}
\usepackage[shortlabels]{enumitem}

\usepackage{amssymb}
\makeatletter
\@ifpackageloaded{amsmath}{}{\usepackage[nosumlimits]{amsmath}}
\makeatother
\usepackage{fixmath} % to provide better bold faced math/greek symbols
\usepackage{xfrac} % for the sfrac command
\usepackage[makeroom]{cancel} % used to mark canceled out terms in equations
\usepackage{dsfont}
\usepackage{mathrsfs}
\usepackage{mathtools}
\usepackage{centernot}

\newcommand{\norm}[2][]{{#1\|}{#2}{#1\|}}
\newcommand{\inprod}[3][]{{#1\langle}{#2},{#3}{#1\rangle}}

\DeclareMathOperator*{\argmin}{arg\,min}

\DeclareMathOperator{\Exp}{\mathbb{E}}

\DeclareMathOperator{\prox}{prox}

\newcommand{\R}{\mathbb{R}}

\usepackage[numbers,square,sort&compress]{natbib}
\usepackage{hyperref} % Must come before cleveref
\usepackage[capitalise]{cleveref} % Use \crefname to define new environment reference. This include needs to come after amsmath include and (of course) before any crefname definition.

\crefformat{equation}{(#2#1#3)}
\crefrangeformat{equation}{(#3#1#4) to~(#5#2#6)}
\crefmultiformat{equation}{(#2#1#3)}{ and~(#2#1#3)}{, (#2#1#3)}{ and~(#2#1#3)}

\crefformat{enumi}{#2#1#3}
\crefrangeformat{enumi}{#3#1#4 to~#5#2#6}
\crefmultiformat{enumi}{#2#1#3}{ and~#2#1#3}{, #2#1#3}{ and~#2#1#3}

\usepackage{amsthm}

\newtheorem{prop}{Proposition}[section]
\numberwithin{prop}{section}

\newtheorem{ass}{Assumption}[section]
\numberwithin{ass}{section}

\newtheorem{lem}{Lemma}[section]
\numberwithin{lem}{section}

\newtheorem{thm}{Theorem}[section]
\numberwithin{thm}{section}

\newtheorem{cor}{Corollary}[section]
\numberwithin{cor}{section}

\newtheorem{defin}{Definition}[section]
\numberwithin{defin}{section}

\numberwithin{rem}{section}

\numberwithin{ex}{section}

\makeatletter
\@ifpackageloaded{cleveref}{
	\crefname{prop}{Proposition}{Propositions}
	\crefname{ass}{Assumption}{Assumptions}
	\crefname{lem}{Lemma}{Lemmas}
	\crefname{thm}{Theorem}{Theorems}
	\crefname{cor}{Corollary}{Corollaries}
	\crefname{defin}{Definition}{Definitions}
	\crefname{rem}{Remark}{Remarks}
	\crefname{ex}{Example}{Examples}
}{}
\makeatother

\graphicspath{{./}{./figures/}}

\usepackage{appendix}
\appendixtitleon % puts \appendixname before each appendix title

\title{Sampling and Update Frequencies in Proximal Variance-Reduced Stochastic Gradient Methods}

\author{
	Martin Morin\\
	Dept. of Automatic Control\\
	Lund University \\
	\url{martin.morin@control.lth.se}
	\And
	Pontus Giselsson\\
	Dept. of Automatic Control\\
	Lund University \\
	\url{pontus.giselsson@control.lth.se}
}

\begin{document}
\maketitle

\begin{abstract}
	Variance-reduced stochastic gradient methods have gained popularity in recent times.
Several variants exist with different strategies for the storing and sampling of gradients and this work concerns the interactions between these two aspects.
We present a general proximal variance-reduced gradient method and analyze it under strong convexity assumptions.
Special cases of the algorithm include SAGA, L-SVRG and their proximal variants.
Our analysis sheds light on epoch-length selection and the need to balance the convergence of the iterates with how often gradients are stored.
The analysis improves on other convergence rates found in the literature and produces a new and faster converging sampling strategy for SAGA.
Problem instances for which the predicted rates are the same as the practical rates are presented together with problems based on real world data.

\end{abstract}

\section{Introduction}
The problem of finding a minimum of a finite sum of functions is common in classification, regression, and general empirical risk minimization.
Each term of the objective is in these cases associated with some error or loss corresponding to a particular data point.
In contemporary problems, the datasets are typically very large and hence the number of terms in the objective function is large.
Traditional iterative minimization algorithms that evaluate the full objective or its gradient each iteration can then become computationally expensive.
Stochastic gradient (SG) methods \cite{robbins_stochastic_1951} have therefore become the methods of choice in this setting \cite{bottou_tradeoffs_2008}, since in each iteration they only evaluate the gradients of a random subset of the terms.

A family of SG methods that have gathered much attention due to their improved convergence properties over ordinary the ordinary SG method are \textit{variance-reduced} SG methods, see \cite{roux_stochastic_2012,johnson_accelerating_2013,defazio_saga_2014,xiao_proximal_2014,schmidt_minimizing_2017,konecny_semi-stochastic_2017,nguyen_sarah_2017,kovalev_dont_2019}.
All variance-reduced methods have a memory over previously evaluated gradients and use them to improve the stochastic estimate of the full gradient.
Although other differences exists, the main separating property between different variance-reduced stochastic gradient method is how the gradient memory is updated.
This work will focus on the effects of how often the memory is updated and of how the stochastic gradient is sampled.

The majority of research into sampling strategies for randomized gradient methods has been on coordinate gradient methods.
Instead of randomly selecting one function from a finite sum, coordinate gradient methods select a random set of coordinates of the gradient and update only those.
One of the first proposed distributions on how these coordinates should be sampled is to sample proportional to a power of the coordinate-wise gradient Lipschitz constant \cite{nesterov_efficiency_2012}.
An arbitrary distribution is allowed in \cite{richtarik_iteration_2014} and \cite{zhao_stochastic_2015} argue that the optimal distribution should be proportional the norm of the coordinate-wise gradient at the current iterate.
Beyond that, \cite{qu_quartz_2015,takac_distributed_2015,csiba_importance_2016,qu_coordinate_2016,richtarik_parallel_2016} present approaches that allow for a combination of randomized mini-batching and arbitrary sampling.

For stochastic gradient and its variance-reduced variants, importance sampling is not as developed.
Variants of importance sampling for the Kaczmarz algorithm and ordinary stochastic gradient are treated in \cite{strohmer_randomized_2008,needell_stochastic_2014}.
For variance-reduced methods, \cite{xiao_proximal_2014} allows for importance sampling in the SVRG setting, while \cite{schmidt_non-uniform_2015} analyzes SAGA under importance sampling.
The results for SAGA are further improved and generalized in \cite{gower_stochastic_2018,qian_saga_2019} to include arbitrary randomized mini-batching with importance sampling.
In this paper, we introduce a general variance-reduced algorithm and prove its linear convergence in the smooth strongly convex regime.
The algorithm allows for importance sampling and have, among others, SAGA \cite{defazio_saga_2014} and L-SVRG \cite{kovalev_dont_2019} as special cases.

The analysis reveals a trade-off between the convergence of the primal iterate (approximate solution) and the dual iterates (stored gradients).
For SAGA, where primal and dual updates are coupled, it is crucial to consider this trade-off when designing samplings and we provide a new sampling strategy that improves on the known convergence rates for SAGA.
For algorithms like L-SVRG, where the memory update is independent of the sampling, it is always beneficial in terms of convergence rate to update more often.
However, this incurs a higher computational cost so we present an update strategy that balances the computational cost against the convergence rate.
Our new rates and computational complexity improve on the previously known results for L-SVRG.

The algorithm in this paper has similarities to the algorithms analyzed in \cite{hofmann_variance_2015} and \cite{zhang_unifying_2019}.
Compared to the memorization algorithm in \cite{hofmann_variance_2015}, our algorithm allows for a proximal term and has a less restrictive gradient memory update.
Our algorithm also allows for importance sampling in SAGA, something that is not supported by the analysis in \cite{zhang_unifying_2019}.
Furthermore, the algorithm of \cite{zhang_unifying_2019} is applied to a larger class of monotone inclusion problems, potentially making the analysis more conservative.

\section{Preliminaries}
Let $\R$ be the set of real numbers.
We will work in finite dimensional real spaces $\R^N$.
Let $\inprod{\cdot}{\cdot}$ denote the standard Euclidean inner product and let $\norm{\cdot}$ be the norm induced by the inner product.
The expected value conditioned on the filtration $\mathcal{F}$ is $\Exp[\cdot | \mathcal{F}]$.
The probability of a discrete random variable taking value $i$ is $P(\cdot = i)$.
We define $\mathbf{1}_{X} = 1$ if the predicate $X$ is true, otherwise $\mathbf{1}_X = 0$.

A convex function $f: \R^d \to \R$ is \emph{$L$-smooth} with $L > 0$ if it is differentiable and its gradient is $\tfrac{1}{L}$-cocoercive, i.e.,
\begin{align*}
	\inprod{\nabla f(x) - \nabla f(y)}{x-y} \geq \tfrac{1}{L} \norm{\nabla f(x) - \nabla f(y)}^2, \quad \forall x,y \in \R^d.
\end{align*}
Note that the definition of smoothness implies \emph{$L$-Lipschitz continuity} of the gradient $\nabla f$.
In fact, for convex $f$, Lipschitz continuity and cocoercivity of $\nabla f$ are equivalent \cite[Corollary 18.17]{bauschke_convex_2017}.
A proper function $f: \R^d \to \R \cup \{\infty\}$ is \emph{$\mu$-strongly convex} with $\mu > 0$ if $f - \frac{\mu}{2}\norm{\cdot}^2$ is convex.

The subdifferential of a $\mu$-strongly convex function is $\mu$-strongly monotone \cite[Example 22.4]{bauschke_convex_2017}, i.e.,
\begin{align*}
	\inprod{u - v}{x-y} \geq \mu \norm{x-y}^2
\end{align*}
holds $\forall x,y \in \mathop{\mathrm{dom}} \partial f$ and $\forall u\in \partial f(x), \forall v \in \partial f(y)$.
A closed, proper and strongly convex function has a unique minimum \cite[Corollary 11.17]{bauschke_convex_2017}.

The \emph{proximal operator} of a closed, convex and proper function $g: \R^d \to \R \cup \{\infty\}$ is defined as
\begin{align*}
	\prox_g(z) = \argmin_x g(x) + \tfrac{1}{2}\norm{x-z}^2.
\end{align*}
Due to strong convexity of $g + \frac{1}{2}\norm{\cdot - z}^2$, the minimum exist and is unique.
Furthermore, the operator $\prox_g$ is \emph{non-expansive}, i.e., Lipschitz continuous with constant $1$ \cite[Proposition 12.28]{bauschke_convex_2017}.

A \emph{Lipschitz distribution} or \emph{Lipschitz sampling} is a probability distribution on $i \in \{1,\dots,n\}$ proportional to the Lipschitz constants $L_i$ of $\nabla f_i$ in \Cref{eq:problem}.

\section{Problem and Algorithm}\label{sec:alg}
We consider the regularized finite sum problem
\begin{align}\label{eq:problem}
	\min_{x \in \R^N} \quad& g(x) + F(x),
\end{align}
where $g: \R^N \rightarrow \R \cup \{\infty\}$ and $F$ is of finite sum form $F(x) = \frac{1}{n} \sum_{i=1}^n f_i(x)$ with $f_i: \R^N \rightarrow \R$ for all $i\in\{1,\dots,n\}$.
We will further make the following assumption on \Cref{eq:problem}.
\begin{ass}[Problem Properties]\label{ass:strong-ass}%
	The function $g: \R^N \rightarrow \R \cup \{\infty\}$ is closed, convex and proper.
	For all $i\in\{1,\dots,n\}$, the function $f_i: \R^N \rightarrow \R$ is convex, differentiable and $L_i$-smooth.
	The function $F: \R^N \rightarrow \R$ is $\mu$-strongly convex, differentiable and $L$-smooth with $L \leq \frac{1}{n}\sum_{i=1}^n L_i$.
\end{ass}
As a consequence of \Cref{ass:strong-ass}, $g + F$ is closed, proper and $\mu$-strongly convex and  hence there exists a unique solution to \Cref{eq:problem}, which we denote $x^\star$.
We propose the following proximal variance-reduced stochastic gradient (PVRSG) method for solving \cref{eq:problem}.
\begin{algorithm}[H]%
	\caption{PVRSG - Proximal Variance-Reduced Stochastic Gradient}\label{alg:pvrsg}
	Given the function $g$, the functions $f_1,\dots,f_n$, initial primal and dual points, $x^0$ and $y_1^0,\dots y_n^0$, iteratively perform the following for $k\in\{0,1,\dots\}$.
	\begin{align*}
		&\textbf{Sampling:}
		&&\text{Randomly sample } (I^k,U_1^k,\dots,U_n^k) \text{ from } \{1,\dots,n\}\times\{0,1\}^n.
		\\
		&\textbf{Primal Update:}
		&&\begin{aligned}[t]
			&\begin{aligned}
				z^{k+1} &=x^k - \tfrac{\lambda}{n} \big(\tfrac{1}{p_{I^k}}(\nabla f_{I^k} (x^k) - y_{I^k}^k) + {\textstyle\sum_{i=1}^n} y_i^k\big), \\
				x^{k+1} &= \prox_{\lambda g}(z^{k+1}).
			\end{aligned}
		\end{aligned}\\
		\\
		&\textbf{Dual Update:}
		&&\begin{aligned}[t]
			&\begin{aligned}
				y_i^{k+1} &= y_i^k + U_i^k(\nabla f_i (x^k) - y_i^k), \quad \forall i \in \{1,\dots,n\}.
			\end{aligned}
		\end{aligned}
	\end{align*}
	The sampling distributions of $(I^k,U_1^k,\dots,U_n^k)$ for all $k\in\{0,1,\dots\}$ are the same and independent.
	Furthermore, the distribution is such that $P(I^k = i) = p_i > 0$ and the expected update frequency $\eta_i > 0$, see \Cref{def:update-freq}, for all $i\in\{1,\dots,n\}$.
	The step-size satisfies $\lambda > 0$.
\end{algorithm}
In \Cref{alg:pvrsg}, the primal variable $x^k$ is updated with a stochastic approximation of the standard proximal gradient (PG) step.
This approximation becomes better the closer the dual variables $y_1^k,\dots,y_n^k$ are to the true gradients $\nabla f_1(x^k),\dots,\nabla f_n(x^k)$.
The purpose of the dual update is to bring these dual variables closer to the true gradients by updating a selection of them with the corresponding gradients at the current iteration.
The more often a dual variable is updated, the closer it will be to the true gradient on average.
We quantify the frequency of the dual updates with the \emph{expected update frequency}, or in short \emph{update frequency}.
\begin{defin}[Expected Update Frequency]\label{def:update-freq}%
	Let $U_1^k,\dots,U_n^k$ be given by the sampling in \Cref{alg:pvrsg}.
	The expected update frequency of the $i$th dual variable is
	\begin{align*}
		\eta_i = \Exp[U_i^k|\mathcal{F}^k],
	\end{align*}
	where $\mathcal{F}^k = \cup_{i=1}^k \mathcal{X}^i$ and $\mathcal{X}^k = \{x^k,y^k,I^{k-1},U_1^{k-1},\dots,U_n^{k-1}\}$.
\end{defin}
Note, the expected update frequency does not depend on the iteration number $k$ since $(I^k,U_1^k,\dots,U_n^k)$ is independently sampled and its distribution does not depend on $k$.

By the nature of the dual update, the dual variables $y_1^k,\dots,y_n^k$ does not necessarily contain gradients evaluated at the same point, i.e., there might not exists $\hat{x}$ such that $\frac{1}{n}\sum_{i=1}^n y_i^k = \nabla F(\hat{x})$.
However, it turns out that if, for all $k\in\{0,1,\dots\}$, there exists such $\hat{x}$, an improved analysis can be made.
This leads to the following assumption.
\begin{ass}[Coherent Dual Update]\label{ass:coherent-dual}%
	For all $i \in \{1,\dots,n\}$, the initial dual variables satisfy $y_i^0 = \nabla f_i(\hat{x})$ for some $\hat{x}$ and it holds that
	\begin{align*}
		U_1^k = U_2^k = \dots = U_n^k
	\end{align*}
	for all $k\in\{0,1,\dots\}$.
\end{ass}

\Cref{alg:pvrsg} contains many special cases with different samplings leading to different algorithms.
The two main algorithms of relevance are SAGA \cite{defazio_saga_2014} and L-SVRG \cite{kovalev_dont_2019}.

\textbf{SAGA:}
SAGA \cite{defazio_saga_2014} only evaluates one gradient each iteration and always save it, i.e., the sampling is defined such that $U_i^k = \mathbf{1}_{i=I^k}$ for all $i\in\{1,\dots,n\}$ which gives the update frequency $\eta_i = p_i$.

\textbf{L-SVRG:}
L-SVRG \cite{kovalev_dont_2019} is inspired by SVRG \cite{johnson_accelerating_2013}, but, instead of a deterministic update of the dual variables, the dual update is based on a weighted coin toss, i.e., $U_i^k = \mathbf{1}_{Q^k < q}$ where $0 < q \leq 1$ and $Q^k$ is independently and uniformly sampled from $[0,1]$.
The expected update frequency is $\eta_i = q$.
\Cref{ass:coherent-dual} is satisfied if the dual variables are initialized in the same point, i.e., there exists $\hat{x}$ s.t. $y_i^0 = \nabla f_i(\hat{x})$ for all $i\in\{1,\dots,n\}$.

We introduce two more special cases to examine the effects of \Cref{ass:coherent-dual} and the expected update frequency $\eta_i$.

\textbf{IL-SVRG (Incoherent Loopless-SVRG):}
IL-SVRG purposefully break the coherent dual assumption, \Cref{ass:coherent-dual}, in L-SVRG.
Each dual variable is independently updated, $U_i^k = \mathbf{1}_{Q_i^k < q}$  where $0 < q \leq 1$ and $Q_i^k$ is independently and uniformly sampled from $[0,1]$.
The update frequency is the same as for L-SVRG, $\eta_i = q$.

\textbf{q-SAGA:}
In q-SAGA \cite{hofmann_variance_2015}, for each iteration, $q \leq n$ indices are sampled uniformly and independently from $\{1,\dots,n\}$ and the corresponding dual variables are updated.
Hence, the sampling is $U_i^k = \mathbf{1}_{i \in J_q}$ where $J_q$ is the set of sampled indices and the update frequency becomes $\eta_i = q/n$.

\section{Convergence Analysis}\label{sec:conv-theory}
We here analyze \Cref{alg:pvrsg} under \Cref{ass:strong-ass} and prove its linear convergence, both with and without the coherent dual update assumption, \Cref{ass:coherent-dual}.
The main results of this analysis can be found in \Cref{thm:pvrsg-conv,thm:pvrsg-conv-coherent} and all proofs will be deferred to \Cref{sec:proof-prop-lemma,sec:proof-theorems}.
Before moving forward with the analysis, we introduce the following necessary quantities that will be used in our Lyapunov analysis.
\begin{defin}\label{def:lyap-terms}%
	Let $x^\star$ be the solution to \Cref{eq:problem} and $y_i^\star = \nabla f_i(x^\star)$ for all $i\in\{1,\dots,n\}$.
	With $y = (y_1,\dots,y_n)$, where $y_i \in \R^n$ for $i \in \{1,\dots,n\}$, we define
	\begin{align*}
		\mathcal{P}(x) = \norm{x - x^\star}^2 - 2\lambda\inprod{\nabla F(x) - \nabla F (x^\star)}{x - x^\star} + \lambda^2\mathcal{V}(x)
	\end{align*}
	and
	\begin{align*}
		\mathcal{D}(y)
		&= \sum_{i=1}^n {\textstyle(1 - \eta_i  + \frac{1}{\gamma_i})\widehat{\gamma}_i} \norm{y_i - y_i^\star}^2
		-{\textstyle (1+\delta^{-1})\lambda^2}\norm{{{\textstyle\frac{1}{n}}\sum_{i=1}^n} y_i -  y_i^\star }^2,
	\end{align*}
	where
	\begin{align*}
		\mathcal{V}(x)
		= \sum_{i=1}^n {\textstyle \frac{(1+\delta)}{n^2 p_i}(\frac{\eta_i\gamma_i}{\delta}  + 1) }\norm{\nabla f_i(x) - \nabla f_i(x^\star)}^2 - \delta\norm{\nabla F(x) - \nabla F(x^\star)}^2
	\end{align*}
	with $\gamma_i \geq 0$, $\delta > 0$ and $\widehat{\gamma}_i = \gamma_i \frac{(1+\delta^{-1}) \lambda^2}{n^2p_i}$.
	If $\gamma_i = 0$ we define $\tfrac{\gamma_i}{\gamma_i} := 1$.
\end{defin}
The variables $\gamma_i$ and $\delta$ are meta-parameters and which will be specified in each proof of the main convergence theorems.
The base of our convergence analysis will be the following proposition.
\begin{prop}\label{prop:main-conv}%
	Let the filtration $\mathcal{F}^k = \cup_{i=1}^k \mathcal{X}^i$ be given by the state $\mathcal{X}^k = \{x^k,y^k,I^{k-1},U_1^{k-1},\dots,U_n^{k-1}\}$.
	If \Cref{ass:strong-ass} holds, the iterates of \Cref{alg:pvrsg} satisfy
	\begin{equation}\label{eq:prop-lyapunov-bound}
		\Exp \big[ \norm{x^{k+1} - x^\star}^2 + \sum_{i=1}^n \widehat{\gamma}_i\norm{y_i^{k+1} - y_i^\star}^2 \big| \mathcal{F}^k \big] \leq  \mathcal{P}(x^k) + \mathcal{D}(y^k),
	\end{equation}
	where $x^\star$ is the unique solution of \Cref{eq:problem} and $y_i^\star = \nabla f_i(x^\star)$.
	See \Cref{def:lyap-terms} for $\mathcal{P}$, $\mathcal{D}$, and $\mathcal{V}$.

	If the primal updates satisfy
	\begin{align}\label{eq:prim-contract}
		\mathcal{P}(x^k) \leq (1-\rho_P)\norm{x^k - x^\star}^2
	\end{align}
	and the dual updates satisfy
	\begin{align}\label{eq:dual-contract}
		\mathcal{D}(y^k) \leq (1-\rho_D)\sum_{i=1}^n \widehat{\gamma}_i\norm{y_i^k - y_i^\star}^2
	\end{align}
	with $\rho_P, \rho_D \in (0,1]$ then \Cref{alg:pvrsg} converges linearly according to
	\begin{align*}
		\Exp \big[ \norm{x^{k} - x^\star}^2 + \sum_{i=1}^n \widehat{\gamma}_i\norm{y_i^{k} - y_i^\star}^2 \big] \in \mathcal{O}((1 - \min(\rho_P, \rho_D))^k).
	\end{align*}
	\begin{proof}
		See \Cref{sec:proof-prop-lemma}.
	\end{proof}
\end{prop}
If we can find algorithm parameters and meta-parameters such that \Cref{alg:pvrsg} satisfy the primal and dual contractions, \Cref{eq:prim-contract} and \Cref{eq:dual-contract} respectively, this proposition proves that \Cref{alg:pvrsg} convergence to a solution.
The following lemma provides these necessary contraction results.
\begin{lem}[Primal Contraction]\label{lem:prim-contract}%
	Let \Cref{ass:strong-ass} hold, the primal iterates of \Cref{alg:pvrsg} satisfy the primal contraction \Cref{eq:prim-contract} with
	\begin{align*}
		\rho_P = \mu\lambda(2 - \nu\lambda)
	\end{align*}
	where $\nu = \max_i {\textstyle  (1+\delta^{-1})\frac{L_i\eta_i\gamma_i}{n p_i}  +  (1+\delta)\frac{L_i}{n p_i} - \delta\mu}$.
	\begin{proof}
		See \Cref{sec:proof-prop-lemma}.
	\end{proof}
\end{lem}
\begin{lem}[Dual Contraction]\label{lem:dual-contract}%
	Let \Cref{ass:strong-ass} and $\gamma_i > 0$ hold for all $i \in\{1,\dots,n\}$, the dual iterates of \Cref{alg:pvrsg} satisfy the dual contraction \Cref{eq:dual-contract} with
	\begin{align*}
		\rho_D = \min_i {\textstyle \eta_i - \tfrac{1}{\gamma_i}}.
	\end{align*}
	\begin{proof}
		See \Cref{sec:proof-prop-lemma}.
	\end{proof}
\end{lem}
\begin{lem}[Dual Contraction - Coherent Updates]\label{lem:dual-contract-uni}%
	Let \Cref{ass:strong-ass} and \ref{ass:coherent-dual} hold and, for all $i\in\{1,\dots,n\}$, let $\gamma_i$ be such that $\gamma_i \geq 0$ and $\frac{L_i}{np_i} \leq \mu$ if $\gamma_i = 0$.
	The dual iterates of \Cref{alg:pvrsg} satisfy the dual contraction \Cref{eq:dual-contract} with
	\begin{align*}
		\rho_D = \min_i
		\begin{cases}
			{\textstyle \eta - (1 - \frac{np_i}{L_i}\mu )\tfrac{1}{\gamma_i}} & \mathrm{if}\text{ } \gamma_i > 0\\
			1 & \mathrm{if}\text{ } \gamma_i = 0
		\end{cases}.
	\end{align*}
	\begin{proof}
		See \Cref{sec:proof-prop-lemma}.
	\end{proof}
\end{lem}
The proofs of our main convergence results, \Cref{thm:pvrsg-conv,thm:pvrsg-conv-coherent}, consist of establishing meta-parameters $\delta$ and $\gamma_i$ for all $i\in\{1,\dots,n\}$ such that these contractions are sufficiently large, i.e., $\min(\rho_P,\rho_D) > 0$.
However, in order to establish these meta-parameters, the following lemma regarding the relationship between the smoothness constants $L_1,\dots,L_n$ and strong convexity constant $\mu$ is needed.
\begin{lem}\label{lem:sampled-lipschitz}%
	Let $L_i$ and $\mu$ be from \Cref{ass:strong-ass} and $p_i$ from \Cref{alg:pvrsg}, then $\max_i \frac{L_i}{np_i} \geq \mu$.
	Furthermore, if $\max_i \frac{L_i}{np_i} = \mu$ then $\frac{L_i}{np_i} = \mu$ for all $i \in \{1,\dots,n\}$.
	\begin{proof}
		See \Cref{sec:proof-prop-lemma}.
	\end{proof}
\end{lem}
The main convergence theorems can now be stated.
\begin{thm}[PVRSG Convergence]\label{thm:pvrsg-conv}%
	Given $\max_i \frac{L_i}{np_i} > \mu$ and \Cref{ass:strong-ass}, if there exists $\rho \in (0, \min_i \eta_i)$ such that
	\begin{align*}
		\rho &= \mu\lambda(2 - \nu\lambda )\\
		\nu &= \min_{\delta > 0} \max_i {\textstyle(1+\delta^{-1})\frac{L_i}{n p_i}\frac{\eta_i}{\eta_i - \rho}  +  (1+\delta)\frac{L_i}{n p_i} - \delta\mu},
	\end{align*}
	then the iterates of \Cref{alg:pvrsg} converge according to
	\begin{align*}
		\Exp \Big[ \norm{x^{k} - x^\star}^2 + \sum_{i=1}^n \widehat{\gamma}_i\norm{y_i^{k} - y_i^\star}^2 \Big] \in \mathcal{O}((1 - \rho)^k)
	\end{align*}
	where $\widehat{\gamma}_1,\dots,\widehat{\gamma}_n$ are given by $\widehat{\gamma}_i ={\textstyle \frac{ \lambda^2}{n^2p_i} \frac{1}{\eta_i - \rho}(1 + \frac{1}{\delta^\star})}$ and $\delta^\star$ is the unique minimizer of the minimization problem defining $\nu$.

	If instead $\max_i \frac{L_i}{np_i} = \mu$ then $\nu = \mu + \mu\max_i \tfrac{\eta_i}{\eta_i-\rho}$ and the convergence is such that
	\begin{align*}
		\Exp \Big[ \norm{x^{k} - x^\star}^2 + \sum_{i=1}^n{\textstyle \frac{ \lambda^2}{n^2p_i} \frac{1}{\eta_i - \tilde{\rho}}} \norm{y_i^{k} - y_i^\star}^2 \Big] \in \mathcal{O}((1 - \tilde{\rho})^k)
	\end{align*}
	holds for all $\tilde{\rho} \in (0, \rho)$.
	\begin{proof}
		See \Cref{sec:proof-theorems}.
	\end{proof}
\end{thm}

\begin{thm}[PVRSG Convergence  - Coherent Dual Updates]\label{thm:pvrsg-conv-coherent}%
	Given \Cref{ass:strong-ass} and \ref{ass:coherent-dual} and $\max_i \frac{L_i}{np_i} > \mu$, if there exists $\rho \in (0, \eta)$ such that
	\begin{align*}
		\rho &= \mu\lambda(2 - \nu \lambda) \\
		\nu &= \mu +  \Big( \max_i {\textstyle \frac{L_i}{n p_i} - \mu\Big)\left(1 + \sqrt{ \frac{\eta}{\eta - \rho}  } \right)^2},
	\end{align*}
	then the iterates of \Cref{alg:pvrsg} converge according to
	\begin{align*}
		\Exp \Big[ \norm{x^{k} - x^\star}^2 + \sum_{i=1}^n \widehat{\gamma}_i\norm{y_i^{k} - y_i^\star}^2 \Big] \in \mathcal{O}((1 - \rho)^k)
	\end{align*}
	where $\widehat{\gamma}_1,\dots,\widehat{\gamma}_n$ are given by $\widehat{\gamma}_i = {\textstyle \frac{ \lambda^2}{n^2p_i}\frac{1}{\eta - \rho} \max(0,1 - \frac{np_i\mu}{L_i})  \Big(1+\sqrt{ \frac{\eta - \rho}{\eta} }\Big) }$.

	If instead $\max_i \frac{L_i}{np_i} = \mu$ then $\nu = \mu$ and $\widehat{\gamma}_i = 0$ for all $i\in\{1,\dots,n\}$. Furthermore, the rate is not restricted to $\rho \in (0,\eta)$, but to $\rho \in (0,1]$.
	\begin{proof}
		See \Cref{sec:proof-theorems}.
	\end{proof}
\end{thm}

Note that the theorems do not provide explicit expressions for the convergence rates, but instead implicitly define them.
Because of this, when we reference the rates of these theorems, we will refer to a numerically computed value.
This computation is done with a combination of convex optimization---for computing $\nu$---and bisection---for finding $\rho$ such that $0 = \rho - \mu\lambda(2 - \nu\lambda)$.

Apart from the coherent dual assumption, our convergence results depend only on the update frequency $\eta_i$ and not on the specifics of the dual sampling that generated it.
Comparing the two theorems, we see that coherent dual updates have greatest effect when the problem is well-conditioned, i.e., when $\tfrac{L_i}{\mu}$ is small for all $i\in\{1,\dots,n\}$.
This stems from the fact that the contraction factor for coherent updates in \Cref{lem:dual-contract-uni} goes towards the contraction factor without coherent updates in \Cref{lem:dual-contract} when $\tfrac{L_i}{\mu}$ increases for all $i\in\{1,\dots,n\}$.

In the extremely well-conditioned case when $\max_i \frac{L_i}{np_i} = \mu$, we see that $\widehat{\gamma}_i = 0$ for all $i\in\{1,\dots,n\}$ and the dual term of the Lyapunov function vanishes completely in \Cref{thm:pvrsg-conv-coherent}.
This is possible due the fact that, in this case, the primal update actually is equal to the true proximal gradient step, regardless of the dual variables.
Also notice that we, as expected, recover the rate for ordinary proximal-gradient.

\section{Special Cases}
In order to provide easily compared rates for SAGA and L-SVRG, we present simplified corollaries of \Cref{thm:pvrsg-conv,thm:pvrsg-conv-coherent} that provide explicit rates.
These rates are by construction conservative compared to the theorems but still improve on previously known best rates.
The corollaries also provide explicit upper bounds on the step-sizes.
Unlike the rates, the bounds are not conservative and match the implicit bounds in \Cref{thm:pvrsg-conv,thm:pvrsg-conv-coherent}.
The proofs of the corollaries are found in \Cref{sec:proof-cor}.

\begin{cor}[SAGA - Conservative Bounds]\label{cor:saga}%
	Given \Cref{ass:strong-ass}, the maximal and recommended step-sizes, $\lambda_{\max}$ and $\lambda^\star$, for SAGA are:
	\begin{align*}
		&\text{If } p_i = \frac{1}{n},
		&&\lambda_{\max} = {\textstyle\frac{2}{C_U L_{\max}}},
		&& \lambda^\star = {\textstyle \frac{2}{C_U L_{\max} + n\mu + \sqrt{ (C_UL_{\max})^2 + (n\mu)^2 }}}. \\[4mm]
		&\text{If } p_i \propto L_i,
		&&\lambda_{\max} = {\textstyle\frac{2}{C_L\bar{L}}},
		&& \lambda^\star = {\textstyle \frac{2}{C_L \bar{L} + p_{\min}^{-1}\mu + \sqrt{ (C_L\bar{L})^2 + (p_{\min}^{-1}\mu)^2 }}}.
	\end{align*}
	where $\bar{L} = \frac{1}{n}\sum_{i=1}^n L_i$, $C_U = {\textstyle 2 + 2\sqrt{1 -\frac{\mu}{L_{\max}}}}$, and $C_L = {\textstyle 2 + 2\sqrt{1 -\frac{\mu}{\bar{L}}}}$.
	The iterates converge with a rate of $\Exp \norm{x^k - x^\star}^2 \in \mathcal{O}\big((1 - \mu\lambda^\star)^k\big)$ when then step-size $\lambda^\star$ is used.
	\begin{proof}
		See \Cref{sec:proof-cor}.
	\end{proof}
\end{cor}

\begin{cor}[L-SVRG - Conservative Bounds]\label{cor:lsvrg}%
	Given \Cref{ass:strong-ass}, the maximal and recommended step-sizes, $\lambda_{\max}$ and $\lambda^\star$, for L-SVRG are:
	\begin{align*}
		&\text{If } p_i = \frac{1}{n},
		&&\lambda_{\max} = {\textstyle\frac{2}{D_U L_{\max}}},
		&& \lambda^\star = {\textstyle \frac{2}{D_U L_{\max} + \eta^{-1}\mu + \sqrt{ (D_UL_{\max})^2 + (\eta^{-1}\mu)^2 }}}. \\[4mm]
		&\text{If } p_i \propto L_i,
		&&\lambda_{\max} = {\textstyle\frac{2}{D_L\bar{L}}},
		&& \lambda^\star = {\textstyle \frac{2}{D_L \bar{L} + \eta^{-1}\mu + \sqrt{ (D_L\bar{L})^2 + (\eta^{-1}\mu)^2 }}}.
	\end{align*}
	where $\bar{L} = \frac{1}{n}\sum_{i=1}^n L_i$, $D_U = 4 - 3 \frac{\mu}{L_{\max}}$ and $D_L = 4 - 3 \frac{\mu}{\bar{L}}$.
	Note that $4 > D_U \geq D_L \geq 1$.
	The iterates converge with a rate of $\Exp \norm{x^k - x^\star}^2 \in \mathcal{O}\big((1 - \mu\lambda^\star)^k\big)$ when then step-size $\lambda^\star$ is used.
	\begin{proof}
		See \Cref{sec:proof-cor}.
	\end{proof}
\end{cor}
The \emph{recommended} step-sizes $\lambda^\star$ are the step-sizes we found that yield the best explicit rates.
However, they are not necessarily optimal w.r.t. the implicit rates in \Cref{thm:pvrsg-conv,thm:pvrsg-conv-coherent}.

% The recommended step-sizes and the resulting rates are all of similar form and they all balance some form of condition number---$\frac{L_{\max}}{\mu}$ or $\frac{\bar{L}}{\mu}$---with the smallest update frequency---$\frac{1}{n}$, $p_{\min}$ and $\eta$.
% We can see this as if they try to balance the primal and dual contractions since the condition number mainly affects primal contraction and the update frequency mainly affects the dual contraction, see \Cref{lem:prim-contract,lem:dual-contract,lem:dual-contract-uni}.
% % However, note that the update frequency can only be chosen independently for L-SVRG and not for SAGA.
% % SAGA therefore runs the risk of being limited by either primal or dual depending on the distribution of Lipschitz constants.

\section{Sampling Design}\label{sec:sampling-design}
Before we present our suggested sampling distributions for SAGA and L-SVRG, we make a few remarks on the parameter selection in \Cref{alg:pvrsg}.

A higher update frequency always yields faster convergence.
However, more frequent dual updates incur a higher computational cost since this require more gradient evaluations.
The  update frequencies therefore needs to be based on the total computational complexity of reaching an $\epsilon$-accurate solution in expectation, i.e., $\Exp \norm{x^k - x^\star}^2 \leq \epsilon$.

The choice of distribution of $p_1,\dots,p_n$ does not change the iteration cost so it can be optimized by only considering the convergence rate, not the computational complexity.
If the update frequencies are uniform, $\eta_i = \eta_j, \forall i,j \in \{1,\dots,n\}$, the meta-parameters $\gamma_1,\dots,\gamma_n$ in \Cref{thm:pvrsg-conv,thm:pvrsg-conv-coherent} also are uniform.
In this case, it can be seen that Lipschitz sampling maximizes the convergence rate, i.e., $p_i \sim L_i$.
However, this is not necessarily true in cases with non-uniform update frequencies.

For SAGA, the expected update frequencies depend on $p_1,\dots,p_n$ and we can therefore not use the optimal choice of uniform update frequencies and Lipschitz sampling of $I^k$.
Instead, we present choice of $p_1,\dots,p_n$ that considers the dependency between the primal and dual update and blends Lipschitz and uniform sampling.
The proposed distribution improves on all other samplings in terms of convergence rate and computational complexity, see \Cref{cor:saga-sampling}.
\begin{cor}[SAGA - Improved Sampling]\label{cor:saga-sampling}%
	Let the sampling distribution and step-size be
	\begin{align*}
		p_i \propto {\textstyle 4L_i + n\mu + \sqrt{(4L_i)^2 + (n\mu)^2}, \quad \lambda = \frac{2}{S}}
	\end{align*}
	where $S = \frac{1}{n}\sum_{i=1}^n (4L_i + n\mu + \sqrt{(4L_i)^2 + (n\mu)^2})$.
	SAGA converges with a rate of $\Exp \norm{x^k - x^\star}^2 \in \mathcal{O}\big( (1 - \mu\lambda)^k \big)$ and achieves an $\epsilon$-accurate solution in expectation within
	\begin{align*}
		\mathcal{O}\big( \tfrac{1}{2}\big({\textstyle \tfrac{1}{n}\sum_{i=1}^n \tfrac{4L_i}{\mu} + n + \sqrt{(\tfrac{4L_i}{\mu})^2 + n^2} }  \big) \log \tfrac{1}{\epsilon}\big)
	\end{align*}
	iterations.
	\begin{proof}
		See \Cref{sec:proof-cor}.
	\end{proof}
\end{cor}

Unlike in SAGA, $\eta_1,\dots,\eta_n$ are always uniform in L-SVRG and can be tuned independently of the primal update.
As remarked on earlier, Lipschitz sampling is then the optimal primal sampling and is therefore used in the following complexity results.
We assume one gradient evaluation is needed in the primal update
\footnote{
	This assumes all dual variables $y_1^k,\dots,y_n^k$ are stored.
	One benefit of PVRSG instances that satisfy \Cref{ass:coherent-dual} is that they can be implemented without storing all dual variables at the cost of one extra gradient evaluation.
	We use the higher memory cost variant in order to compare to SAGA under equal memory requirements.
}
and that, in expectation, $n\eta$ are needed in the dual update.
\begin{cor}[L-SVRG - Computational Complexity]\label{cor:lsvrg-complexity-full}%
	Let Lipschitz sampling---$p_i \sim L_i$ for all $i\in\{1,\dots,n\}$---and the step-size from \Cref{cor:lsvrg} be used.
	L-SVRG achieves an $\epsilon$-accurate solution within
	\begin{align*}
		\mathcal{O}\Big( (1 + n\eta)\Big(D_L \tfrac{\bar{L}}{\mu} + \tfrac{1}{\eta}\Big) \log \tfrac{1}{\epsilon} \Big)
	\end{align*}
	iterations where $\bar{L} = \tfrac{1}{n}\sum_{i=1}^n L_i$ and $D_L$ is given by \Cref{cor:lsvrg}.
	The expected update frequency that minimizes the complexity, and the corresponding complexity, are
	\begin{align*}
		\eta^\star = \sqrt{\tfrac{\mu}{nD_L\bar{L}}}
		\quad \text{and} \quad
		\mathcal{O}\Big(\Big(\sqrt{n} + \sqrt{D_L \tfrac{\bar{L}}{\mu}}\Big)^2 \log \tfrac{1}{\epsilon}\Big).
	\end{align*}
	\begin{proof}
		See \Cref{sec:proof-cor}.
	\end{proof}
\end{cor}
The complexity of L-SVRG in $\Cref{cor:lsvrg-complexity-full}$ is worse than that of SAGA in \Cref{cor:saga-sampling} when $n > 2$.
The cheaper iteration cost of SAGA clearly outweighs loss of the coherent dual update, \Cref{ass:coherent-dual}.
With the choice of update frequency for L-SVRG in \cref{cor:lsvrg-complexity-full}, the expected time between dual updates is $\frac{1}{\eta^\star} \propto \sqrt{n \tfrac{\bar{L}}{\mu}}$.
This is in contrast to most results for SVRG and L-SVRG that have epoch lengths proportional to either $n$ or $\tfrac{L}{\mu}$ \cite{johnson_accelerating_2013,xiao_proximal_2014,babanezhadharikandehStopWastingMyGradients2015,sebbouhClosingGapTheory2019,kovalev_dont_2019}.

\section{Numerical Experiments}\label{sec:experiments}

All algorithms have been implemented in \texttt{Julia} \cite{bezansonJuliaFreshApproach2017} and can be found at \url{https://github.com/mvmorin/VarianceReducedSG.jl}.

\begin{figure}[t]
\begin{minipage}[t]{.49\textwidth}
\centering
\includegraphics[width=\linewidth]{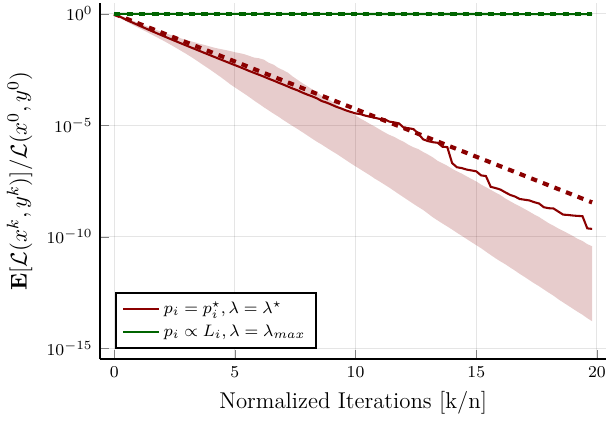}
{(a) SAGA}
\end{minipage}
\begin{minipage}[t]{.49\textwidth}
\centering
\includegraphics[width=\linewidth]{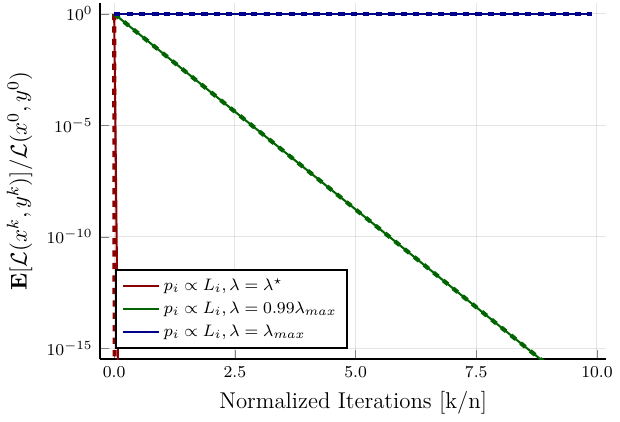}
{(b) L-SVRG}
\end{minipage}\\[.8em]
\begin{minipage}[t]{.49\textwidth}
\centering
\includegraphics[width=\linewidth]{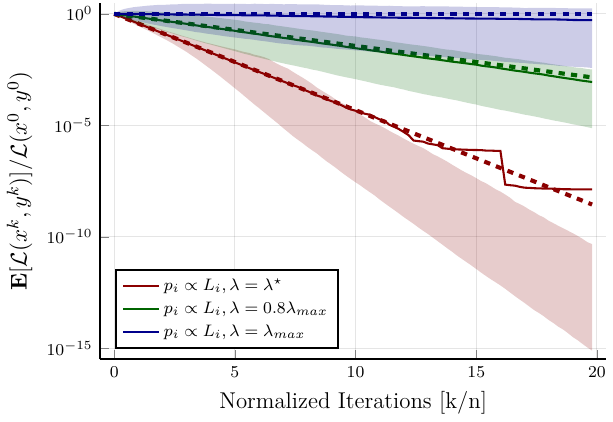}
{(c) q-SAGA}
\end{minipage}
\begin{minipage}[t]{.49\textwidth}
\centering
\includegraphics[width=\linewidth]{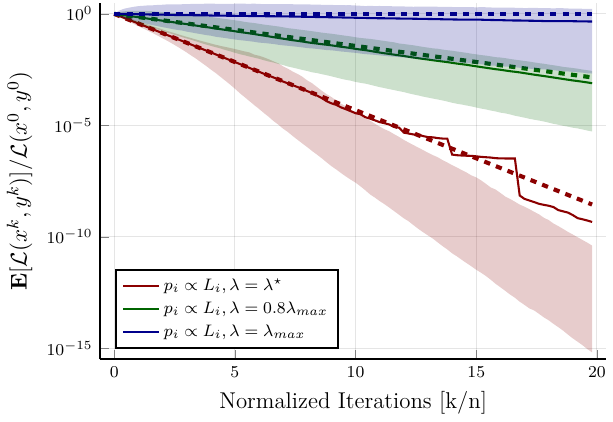}
{(d) IL-SVRG}
\end{minipage}\\
\caption[]{
One-dimensional least squares.
The expected value $\Exp \mathcal{L}(x^k,y^k)$ is estimated with the average of 10000 runs where $\mathcal{L}(x,y)$ is taken from \Cref{thm:pvrsg-conv,thm:pvrsg-conv-coherent}.
The shaded areas represent the 5-95 percentile of the runs.
The dashed lines are the predicted rates.
The step-sizes $\lambda^\star$ and $\lambda^{\max}$ are the optimal and maximal step-sizes according to \Cref{thm:pvrsg-conv,thm:pvrsg-conv-coherent}.
The expected update frequency is $\eta = \frac{1}{n}$ for all algorithms except SAGA where it depends on the sampling $p_i$.
The sampling $p_i^\star$ is from \Cref{cor:saga-sampling}.
%\footnotemark
}
\label{fig:tight_saga_lsvrg_qsaga_ilsvrg}
\end{figure}

\textbf{Simple Least Squares}
The analysis predicts performance accurately for a one-dimensional least squares problem,
\begin{align*}
\min_x \tfrac{1}{n} \sum_{i=1}^n (a_ix - b_i)^2.
\end{align*}
A comparison of theoretical and practical rates for this problem is found in \Cref{fig:tight_saga_lsvrg_qsaga_ilsvrg}.
The data $a_i$ and $b_i$ have been independently drawn from a unit normal distribution and the number of functions is $n = 100$.

For L-SVRG, \Cref{fig:tight_saga_lsvrg_qsaga_ilsvrg} shows fast convergence and very narrow 5-95 percentile---it is not even visible.
This is due to the $\max_i \frac{L_i}{np_i} = \mu$ condition being satisfied and then the gradient estimate is exact.
Since the condition number of the problem is equal to $1$, it is possible to solve the problem in one iteration.

For SAGA, we see in \Cref{fig:tight_saga_lsvrg_qsaga_ilsvrg} that both the maximal and optimal step-sizes are predicted well.
However, note that the sampling distribution $p_1,\dots,p_n$ are not the same for the two cases.

Comparing q-SAGA and IL-SVRG in \Cref{fig:tight_saga_lsvrg_qsaga_ilsvrg}, we see similar performance.
This was predicted by \Cref{thm:pvrsg-conv} since, despite the dual updates being different, the algorithms have the same expected update frequency.
Comparing to L-SVRG in \Cref{fig:tight_saga_lsvrg_qsaga_ilsvrg} we see the huge impact of the coherent dual assumption in this very well-conditioned case.

\begin{figure}[t]
\begin{minipage}[t]{.49\textwidth}
\centering
\includegraphics[width=\linewidth]{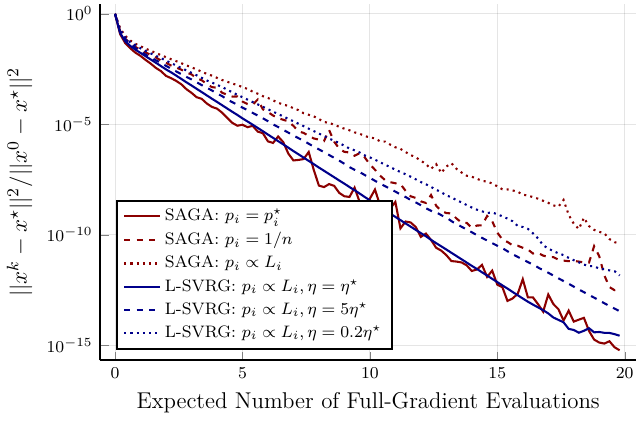}
{(a) Lasso: \texttt{ijcnn1}, $\xi=0.001$,\\ $n=49990$, $\tfrac{\bar{L}}{\mu} \approx 165 $}
\end{minipage}
\begin{minipage}[t]{.49\textwidth}
\centering
\includegraphics[width=\linewidth]{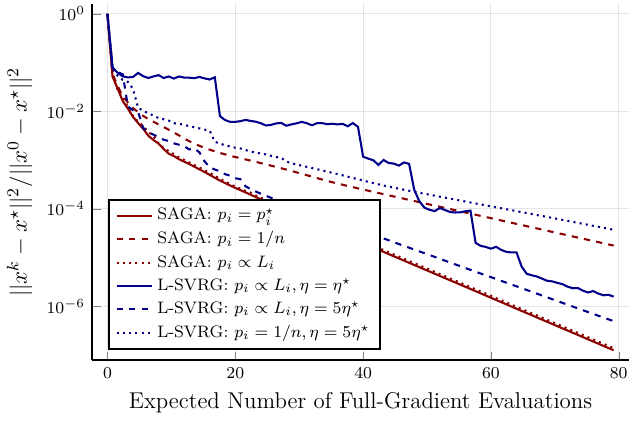}
{(b) Lasso: \texttt{protein}${}^\ast$, $\xi=0.0001$,\\ $n=17766$, $\tfrac{\bar{L}}{\mu} \approx 3\cdot 10^{6}$}
\end{minipage}\\[.8em]
\begin{minipage}[t]{.49\textwidth}
\centering
\includegraphics[width=\linewidth]{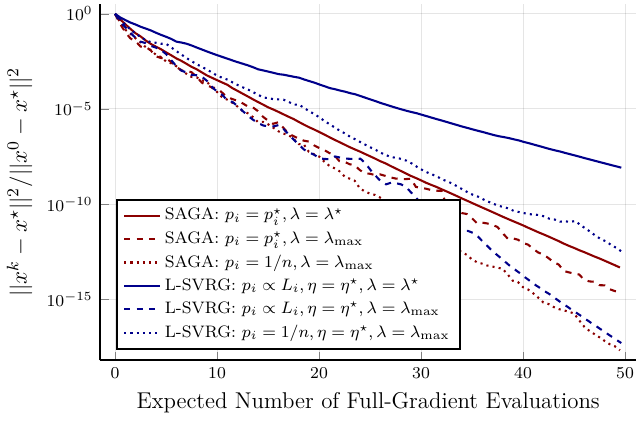}
{(c) Least Squares: \texttt{splice}, $\xi=0$,\\ $n=1000$, $\tfrac{\bar{L}}{\mu} \approx 806$}
\end{minipage}
\begin{minipage}[t]{.49\textwidth}
\centering
\includegraphics[width=\linewidth]{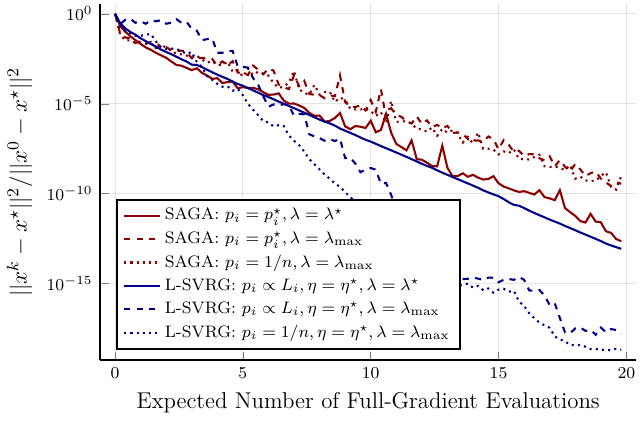}
{(d) Least Squares: \texttt{acoustic}, $\xi=0$,\\ $n=78823$, $\tfrac{\bar{L}}{\mu} \approx 1480$}
\end{minipage}\\
\caption[]{
The \emph{Expected Number of Full-Gradient Evaluations} is the number of full gradient evaluations the algorithms is expected to perform in $k$ iterations, $\frac{1}{n}k$ for SAGA and $(\frac{1}{n} + \eta)k$ for L-SVRG.
The step-sizes, $\lambda^\star$ or $\lambda_{\max}$ are taken from the corresponding result in \Cref{cor:saga}-\ref{cor:lsvrg}, $\lambda^\star$ is used if no step-size is given.
The sampling $p_i^\star$ is from \Cref{cor:saga-sampling} and $\eta^\star$ is the update frequency from \Cref{cor:lsvrg-complexity-full}.
The average condition number $\bar{L} = \frac{1}{n}\sum_{i=1}^n L_i$.
The \texttt{protein} dataset has a feature consisting of only zeros.
\texttt{protein}$^\ast$ has this feature removed in order to preserve strong convexity.
%\footnotemark
}
\label{fig:libsvm-experiments}
\end{figure}
%\footnotetext{For color figures we refer to the online version of the paper.}

\textbf{Lasso Problem}
Here we consider a Lasso regression problem of the form
\begin{align*}
\min \norm{Ax-b}_2^2 + \xi \norm{x}_1,
\end{align*}
where the matrix $A$ and vector $b$ consist of the features and classes from different datasets from the \texttt{LibSVM} database \cite{chang_libsvm_2011}.
The regularization parameter $\xi$ is tuned for each problem such that the solution have roughly 15-20\% sparsity.
SAGA and L-SVRG with different sampling and update frequencies are compared in \Cref{fig:libsvm-experiments}(a)-(b).
L-SVRG was tested with $\eta \in \{0.2\eta^\star, \eta^\star, 5\eta^\star\}$ and either Lipschitz or uniform sampling and the three best perform configuration are shown in \Cref{fig:libsvm-experiments}.
Further comparisons between SAGA and L-SVRG with larger step-size choices can be found in \Cref{fig:libsvm-experiments}(c)-(d).
In these experiments the regularization parameter was set to $\xi = 0$.

We see that it is sometime possible achieve better convergence rate by deviating from the optimal parameter choices in \Cref{cor:saga,cor:lsvrg,cor:saga-sampling,cor:lsvrg-complexity-full}.
However, these are single realizations of random processes and there will be variance between runs, especially for the larger step-sizes.
The consistency of our suggested parameter choices should be noted though.
Especially SAGA with the sampling from \Cref{cor:saga-sampling} are always among the better alternatives.

Note the slow step-like convergence of L-SVRG with $\eta = \eta^\star$ in the \texttt{protein}$^\ast$ example.
The faster convergence of $\eta = 5\eta^\star$ suggest that $\eta^\star$ does not properly balance the primal and dual updates.
Since we perform a worst case analysis, there are many reasons for why this might be the case.
Our analysis also only focus on asymptotic linear rates and does not capture transient behavior.
The last point is especially important when considering very ill-conditioned or maybe even non-strongly convex problems.
In these cases, the transient phase is the most important part since the achievable linear rates are very small or zero.

\section{Conclusion}\label{sec:conclusion}
A general stochastic variance-reduced gradient method has been analyzed and problems have been presented where the predicted rates are close to real world rates.
We have demonstrated the need to balance the updates of the primal and dual variables.
For L-SVRG, we presented a new condition number dependent update probability for the dual variables.
For SAGA, and other methods where the dual update depends on the primal update, the primal sampling needs to consider both updates.
Lipschitz sampling, which appears to be optimal for methods with independent dual updates, can for SAGA lead to slow convergence.
We have presented a new sampling for SAGA that balances the primal and dual update and consistently performs well.

\paragraph{\normalfont{\textbf{Funding}}}
This work is funded by the Swedish Research Council via grant number 2016-04646.

\begin{appendices}
	\section{Proofs of Proposition and Lemmas}\label{sec:proof-prop-lemma}

\begin{proof}[\textbf{Proof of \Cref{prop:main-conv}}]%
	Let $x^\star$ be the unique solution to \Cref{eq:problem}.
	With $g$ being proper, closed and convex, the primal updates satisfy
	\begin{equation}\label{eq:prim-update-lyap}
		\begin{aligned}
			&\Exp[\norm{x^{k+1} - x^\star}^2|\mathcal{F}^k] \\
			&\quad= \Exp[\norm{\prox_{\lambda g}(z^{k+1}) - \prox_{\lambda g}(x^\star - \lambda \nabla F(x^\star))}^2|\mathcal{F}^k] \\
			&\quad\leq
			%\Exp[\norm{z^{k+1} - (x^\star - \lambda \nabla F(x^\star))}^2|\mathcal{F}^k] \\
			%&\quad=
			\Exp[\norm{x^k - \tfrac{\lambda}{n}\big(\tfrac{1}{p_{I^k}}(\nabla f_{I^k}(x^k) - y_{I^k}^k) + \textstyle\sum_{i=1}^n y_i^k\big) - x^\star + \lambda \nabla F(x^\star)}^2|\mathcal{F}^k] \\
			&\quad= \norm{(x^k - \lambda \nabla F(x^k)) - (x^\star - \lambda \nabla F(x^\star))}^2\\
		&\quad\quad + \lambda^2 \Exp\big[\norm[\big]{\big(\tfrac{1}{n p_{I^k}} \nabla f_{I^k}(x^k) - \nabla F(x^k)\big) - \big(\tfrac{1}{n p_{I^k}}y_{I^k}^k - {\textstyle\frac{1}{n}\sum_{i=1}^n} y_i^k \big)}^2|\mathcal{F}^k\big] \\
		&\quad= \norm{x^k - x^\star}^2 - 2\lambda\inprod{\nabla F(x^k) - \nabla F(x^\star)}{x^k-x^\star} + \lambda^2\norm{\nabla F(x^k) -\nabla F(x^\star)}^2\\
	&\quad\quad + \lambda^2 \Exp\big[\norm[\big]{\big(\tfrac{1}{n p_{I^k}} \nabla f_{I^k}(x^k) - \nabla F(x^k)\big) - \big(\tfrac{1}{n p_{I^k}}y_{I^k}^k - {\textstyle\frac{1}{n}\sum_{i=1}^n} y_i^k \big)}^2|\mathcal{F}^k\big].
\end{aligned}
\end{equation}
The first equality is given by the solution being a fixed point to the proximal-gradient update $x^\star = \prox_{\lambda g}(x^\star - \lambda \nabla F(x^\star))$ \cite[Corollary 28.9]{bauschke_convex_2017}.
The first inequality is due to the non-expansiveness of $\prox_{\lambda g}$.
The second to last equality is given by $\Exp\norm{X}^2 =\norm{\Exp X}^2 + \Exp\norm{X - \Exp X}^2$ where $X$ is a random variable.

The last term in \Cref{eq:prim-update-lyap} satisfies the following upper bound for all $\delta > 0$:
\begin{equation}\label{eq:update-variance-bound}
	\begin{aligned}
		&\Exp[\norm{(\tfrac{1}{n p_{I^k}} \nabla f_{I^k}(x^k) - \nabla F(x^k)) - (\tfrac{1}{n p_{I^k}}y_{I^k}^k - {\textstyle\frac{1}{n}\sum_{i=1}^n} y_i^k )}^2|\mathcal{F}^k]\\
		&\quad\begin{aligned}
			= \Exp[\| &\tfrac{1}{np_{I^k}}(\nabla f_{I^k}(x^k) - \nabla f_{I^k}(x^\star)) - (\nabla F(x^k) - \nabla F(x^\star))\\
			&- \tfrac{1}{np_{I^k}}(y_{I^k}^k - y_{I^k}^\star) + {\textstyle\frac{1}{n}\sum_{i=1}^n}(y_i^k - y_i^\star) \|^2 |\mathcal{F}^k]
		\end{aligned} \\
		&\quad\leq (1+\delta)\Exp[\norm{ \tfrac{1}{np_{I^k}}(\nabla f_{I^k}(x^k) - \nabla f_{I^k}(x^\star)) - (\nabla F(x^k) - \nabla F(x^\star)) }^2 |\mathcal{F}^k]\\
		&\quad\quad + (1+\delta^{-1})\Exp[\norm{\tfrac{1}{np_{I^k}}(y_{I^k}^k - y_{I^k}^\star) - {\textstyle\frac{1}{n}\sum_{i=1}^n}(y_i^k - y_i^\star)}^2 |\mathcal{F}^k]\\
		&\quad= (1+\delta)\big(\Exp[\norm{\tfrac{1}{np_{I^k}}(\nabla f_{I^k}(x^k) - \nabla f_{I^k}(x^\star))}^2|\mathcal{F}^k] - \norm{\nabla F(x^k) - \nabla F(x^\star)}^2\big) \\
		&\quad\quad + (1+\delta^{-1})\big(\Exp[\norm{\tfrac{1}{np_{I^k}}(y_{I^k}^k - y_{I^k}^\star)}^2|\mathcal{F}^k] - \norm{{\textstyle\frac{1}{n}\sum_{i=1}^n}(y_i^k - y_i^\star)}^2\big) \\
		&\quad= (1+\delta){\textstyle\sum\frac{1}{n^2 p_i}}\norm{\nabla f_i(x^k) - \nabla f_i(x^\star)}^2 - (1+\delta)\norm{\nabla F(x^k) - \nabla F(x^\star)}^2 \\
		&\quad\quad + (1+\delta^{-1}){\textstyle\sum\frac{1}{n^2 p_i}}\norm{y_i^k - y_i^\star}^2 - (1+\delta^{-1})\norm{{\textstyle\frac{1}{n}\sum_{i=1}^n}(y_i^k - y_i^\star)}^2.
	\end{aligned}
\end{equation}
The inequality is given by Young's inequality, the second to last equality is given by $\Exp\norm{X - \Exp X}^2 = \Exp\norm{X}^2-\norm{\Exp X}^2$ where $X$ is a random variable.

The dual updates satisfy
\begin{align*}
	\norm{y_i^{k+1} - y_i^\star}^2
	&= \norm{y_i^k + U_i^k(\nabla f_i(x^k) - y_i^k) - y_i^\star}^2 \\
	&= (1- U_i^k)\norm{y_i^k - y_i^\star}^2 + U_i^k\norm{\nabla f_i(x^k) - \nabla f_i(x^\star)}^2
\end{align*}
since $U_i^k \in \{0,1\}$.
Summing over all terms, taking expected value and using linearity of the expected value give
\begin{equation}\label{eq:dual-update-lyap}
	\begin{aligned}
		\Exp[{\textstyle\sum_{i=1}^n} \widehat{\gamma}_i\norm{y_i^{k+1} - y_i^\star}^2|\mathcal{F}^k]
		&= {\textstyle\sum_{i=1}^n} (1-\eta_i)\widehat{\gamma}_i\norm{y^k - y^\star}^2 \\
		&\quad + {\textstyle\sum_{i=1}^n}{\eta_i \widehat{\gamma}_i}\norm{\nabla f_i(x^k) - \nabla f_i(x^\star)}^2.
	\end{aligned}
\end{equation}
Adding \Cref{eq:dual-update-lyap} to \Cref{eq:prim-update-lyap} and substituting in \Cref{eq:update-variance-bound} and using the definition $\frac{\gamma_i}{\gamma_i} = 1$ when $\gamma_i = 0$ then yield \Cref{eq:prop-lyapunov-bound}.
Applying \Cref{eq:prim-contract} and \Cref{eq:dual-contract}, using the law of total expectation and telescoping the inequalities give the stated rate.
\end{proof}

\begin{proof}[\textbf{Proof of \Cref{lem:prim-contract}}]%
	First note that $\mu$-strong monotonicity and the Cauchy--Schwarz inequality imply
	\begin{equation}\label{eq:reverse-lipschitz}
		\norm{\nabla F(x^k) - \nabla F(x^\star)} \geq \mu \norm{x^k - x^\star}.
	\end{equation}
	Consider the terms of $\mathcal{V}(x^k)$.
	Using \Cref{eq:reverse-lipschitz} and Cauchy--Schwarz in the last term yield
	\begin{align*}
		\norm{\nabla F(x^k) - \nabla F(x^\star)}^2
		&\geq \mu\norm{\nabla F(x^k) - \nabla F(x^\star)}\norm{x^k - x^\star} \\
		&\geq \mu\inprod{\nabla F(x^k) - \nabla F(x^\star)}{x^k - x^\star}.
	\end{align*}
	Using $\frac{1}{L_i}$-cocoercivity of $\nabla f_i$ in the first term of $\mathcal{V}(x^k)$ yields
	\begin{align*}
		\sum_{i=1}^n & {\textstyle (\frac{\eta_i\gamma_i}{\delta}  + 1)\frac{1}{n^2 p_i} }\norm{\nabla f_i(x^k) - \nabla f_i(x^\star)}^2 \\
		&\leq \sum_{i=1}^n {\textstyle (\frac{\eta_i\gamma_i}{\delta}  + 1)\frac{L_i}{n^2 p_i} } \inprod{\nabla f_i(x^k) - \nabla f_i(x^\star)}{x^k - x^\star} \\
		&\leq \max_i({\textstyle(\frac{\eta_i\gamma_i}{\delta}  + 1)\frac{L_i}{n p_i}}) \inprod{\nabla F(x^k) - \nabla F(x^\star)}{x^k - x^\star}.
	\end{align*}
	Adding the terms back together yields
	\begin{align*}
		\mathcal{V}(x^k)
		&\leq \lambda^2 \nu \inprod{\nabla F(x^k) - \nabla F(x^\star)}{x^k - x^\star}
	\end{align*}
	where $\nu = \max_i(  (1+\delta^{-1})\frac{L_i\eta_i\gamma_i}{n p_i}  +  (1+\delta)\frac{L_i}{n p_i} - \delta\mu)$.
	This can now be summarized as
	\begin{align*}
		\mathcal{P}(x^k)
		&\leq \norm{x^k - x^\star}^2  - \lambda\big(2 - \nu\lambda \big) \inprod{\nabla F(x^k) - \nabla F (x^\star)}{x^k - x^\star} \\
		&\leq \norm{x^k - x^\star}^2  - \mu\lambda\big(2 - \nu\lambda \big) \norm{x^k - x^\star}^2 \\
		&= (1 - \rho_P)\norm{x^k - x^\star}^2,
	\end{align*}
	where $\rho_P = \mu\lambda(2 - \nu\lambda)$ and the last inequality is given by the strong monotonicity of $\nabla F$.
\end{proof}

\begin{proof}[\textbf{Proof of \Cref{lem:dual-contract}}]%
	Since norms are non-negative, we have
	\begin{align*}
		\mathcal{D}(y^k)
		\leq \sum_{i=1}^n {\textstyle(1 - \eta_i  + \frac{1}{\gamma_i})\widehat{\gamma}_i} \norm{y_i^k - y_i^\star}^2
		\leq (1 - \rho_D)\sum_{i=1}^n \widehat{\gamma}_i \norm{y_i^k - y_i^\star}^2,
	\end{align*}
	where $\rho_D = \min_i(\eta_i - \tfrac{1}{\gamma_i})$.
\end{proof}

\begin{proof}[\textbf{Proof of \Cref{lem:dual-contract-uni}}]%
	From \Cref{ass:coherent-dual}, we know that there exists $\phi^k$ such that $y_i^k = \nabla f_i(\phi^k) ,\forall i \in \{1,\dots,n\}$.
	Using this and $y_i^\star = \nabla f_i(x^\star)$ yield
	\begin{gather*}
		\norm{{{\textstyle\frac{1}{n}}\sum_{i=1}^n}\big[y_i^k  -  y_i^\star\big]}^2
		= \norm{\nabla F(\phi^k) - \nabla F(x^\star)}^2
		\geq \mu\norm{\nabla F(\phi^k) - \nabla F(x^\star)}\norm{\phi^k-x^\star} \\
		\geq \mu\inprod{\nabla F(\phi^k) - \nabla F(x^\star)}{\phi^k-x^\star}
		= \mu\frac{1}{n}\sum_{i=1}^n \inprod{\nabla f_i(\phi^k) - \nabla f_i(x^\star)}{\phi^k-x^\star} \\
		\geq \mu\frac{1}{n}\sum_{i=1}^n \frac{1}{L_i}\norm{\nabla f_i(\phi^k) - \nabla f_i(x^\star)}^2
		= \sum_{1=1}^n \frac{\mu}{nL_i} \norm{y_i^k - y_i^\star}^2,
	\end{gather*}
	where the inequalities are given by $\mu$-strong monotonicity of $\nabla F$, Cauchy--Schwarz, and $\frac{1}{L_i}$-cocoercivity of $\nabla f_i$.
	Inserting this into $\mathcal{D}(y^k)$ and using that $\frac{\widehat{y}_i}{y_i} = \frac{(1+\delta^{-1})\lambda^2}{n^2p_i}$ for all $i\in\{1,\dots,n\}$---note that we defined $\frac{\gamma_i}{\gamma_i} = 1$ if $\gamma_i = 0$---
	give
	\begin{align*}
		\mathcal{D}(y^k)
		\leq \sum_{i=1}^n {\textstyle(1 - \eta  + (1 - np_i\frac{\mu}{L_i})\frac{1}{\gamma_i})\widehat{\gamma}_i} \norm{y_i^k - y_i^\star}^2.
	\end{align*}
	For each term we see that if $\gamma_i > 0$ then
	\begin{align*}
		{\textstyle(1 - \eta  + (1 - np_i\frac{\mu}{L_i})\frac{1}{\gamma_i})\widehat{\gamma}_i} \norm{y_i^k - y_i^\star}^2
		\leq (1 - \rho_i^+) \widehat{\gamma}_i \norm{y_i^k - y_i^\star}^2
	\end{align*}
	with $\rho_i^+ = \eta  - (1 - np_i\frac{\mu}{L_i})\frac{1}{\gamma_i}$.
	If $\gamma_i = 0$, then $\frac{L_i}{np_i} \leq \mu$ and
	\begin{align*}
		&{\textstyle(1 - \eta  + (1 - np_i\frac{\mu}{L_i})\frac{1}{\gamma_i})\widehat{\gamma}_i} \norm{y_i^k - y_i^\star}^2 \\
		&\quad=
		{\textstyle(1 - np_i\frac{\mu}{L_i}) \frac{(1+\delta^{-1})\lambda^2}{n^2p_i}} \norm{y_i^k - y_i^\star}^2
		\leq 0
		\leq (1-1)\widehat{\gamma}_i \norm{y_i^k - y_i^\star}^2.
	\end{align*}
	This gives $\mathcal{D}(y^k) \leq (1-\rho_D)\sum_{i=1}^n \widehat{\gamma}_i \norm{y_i^k - y_i^\star}^2$ where
	\begin{align*}
		\rho_D =
		\begin{cases}
			{\textstyle \eta - (1 - \frac{np_i}{L_i}\mu )\tfrac{1}{\gamma_i}} & \mathrm{if}\text{ } \gamma_i > 0\\
			1 & \mathrm{if}\text{ } \gamma_i = 0 \\
		\end{cases}.
	\end{align*}
\end{proof}

\begin{proof}[\textbf{Proof of \Cref{lem:sampled-lipschitz}}]%
	An $L$-smooth and $\mu$-strongly convex function must satisfy $L \geq \mu$.
	Assuming $\max_i(\frac{L_i}{np_i}) < \mu$ yields the following contradiction
	\begin{align*}
		\mu > \max_i \tfrac{L_i}{np_i} = \sum_{j=1}^n  p_j \max_i \tfrac{L_i}{np_i} \geq \sum_{j=1}^n  p_j \tfrac{L_j}{np_j} = \sum_{i=1}^n \tfrac{L_i}{n} \geq L.
	\end{align*}
	If $\max_i(\frac{L_i}{np_i}) = \mu$, equality must hold everywhere and we have
	\begin{align*}
		0 = \sum_{j=1}^n  p_j \max_i \tfrac{L_i}{np_i} - \sum_{j=1}^n  p_j \tfrac{L_j}{np_j} = \sum_{j=1}^n  p_j (\max_i \tfrac{L_i}{np_i} - \tfrac{L_j}{np_j}).
	\end{align*}
	Since $p_j > 0$ and $\max_i \tfrac{L_i}{np_i} - \tfrac{L_j}{np_j} \geq 0$, we have $\max_i \tfrac{L_i}{np_i} = \tfrac{L_j}{np_j}$ for all $j\in\{1,\dots,n\}$.
\end{proof}

\section{Proofs of Theorems}\label{sec:proof-theorems}

\begin{proof}[\textbf{Proof of \Cref{thm:pvrsg-conv}}]%
	Application of \Cref{lem:prim-contract} and \ref{lem:dual-contract} in \Cref{prop:main-conv} yields the convergence rate
	\begin{align*}
		\Exp \Big[ \norm{x^{k} - x^\star}^2 + \sum_{i=1}^n \widehat{\gamma}_i\norm{y_i^{k} - y_i^\star}^2 \Big] \in \mathcal{O}((1 - \min(\rho_P, \rho_D))^k)
	\end{align*}
	with
	\begin{align*}
		\rho_P &= \mu\lambda(2- \nu\lambda) \\
		\rho_D &= \min_i \eta_i - \tfrac{1}{\gamma_i} \\
		\nu &= \max_i  {\textstyle(1+\delta^{-1})\frac{L_i\eta_i\gamma_i}{n p_i}  +  (1+\delta)\frac{L_i}{n p_i} - \delta\mu},
	\end{align*}
	which hold for all choices of $\delta > 0$ and $\gamma_i > 0$ for all $i\in\{1,\dots,n\}$.
	If there exists $\delta$ and $\gamma_1,\dots,\gamma_n$ such that $\min(\rho_P,\rho_D) \in (0,1]$ we have convergence.
	We restrict ourselves to only search for $\delta$ and $\gamma_1,\dots,\gamma_n$ such that $\rho_P = \rho_D = \rho$ for some $\rho \in (0,1]$.
	For all $i\in\{1,\dots,n\}$, select $\gamma_i = \frac{1}{\eta_i - \rho}$, which is positive when $\rho < \eta_i$, and convergence is then proved if there exists $\rho \in (0,\min_i \eta_i)$ and $\delta > 0$ such that
	\begin{align*}
		\rho &= \mu\lambda(2 - \nu\lambda )\\
		\nu &= \max_i {\textstyle(1+\delta^{-1})\frac{L_i}{n p_i}\frac{\eta_i}{\eta_i - \rho}  +  (1+\delta)\frac{L_i}{n p_i} - \delta\mu}.
	\end{align*}
	The variable $\nu$ can be minimized w.r.t. $\delta$ if $\max_i \frac{L_i}{np_i} > \mu$.
	The minimum then exists and is unique since $\nu$ as a function of $\delta$ is continuous, strictly convex, and $\nu\to\infty$ both when $\delta \to 0^+$ and $\delta \to \infty$.
	Calling the minimum point $\delta^\star$, noting that $\delta^\star > 0$, and inserting it and the choice of $\gamma_i$ in the expression for $\widehat{\gamma}_i$ from \Cref{def:lyap-terms} yield the first statement of the theorem.

	When $\max_i \frac{L_i}{np_i} \not> \mu$, \Cref{lem:sampled-lipschitz} gives $\frac{L_i}{n p_i} = \mu$ for all $i\in\{1,\dots,n\}$ and
	\begin{align*}
		\nu = \mu + \mu(1+\delta^{-1})\max_i \tfrac{\eta_i}{\eta_i-\rho}.
	\end{align*}
	This can not be minimized w.r.t. $\delta$ since the $\inf$ is not attained.
	However, any $\delta>0$ will yield a valid $\rho$ and $\widehat{\gamma}_i$, giving the rate
	\begin{align*}
		&\Exp \Big[ \norm{x^{k} - x^\star}^2 + \sum_{i=1}^n{\textstyle \frac{ \lambda^2}{n^2p_i} \frac{1}{\eta_i - \rho}} \norm{y_i^{k} - y_i^\star}^2 \Big] \\
		&\quad \leq \Exp \Big[ \norm{x^{k} - x^\star}^2 + \sum_{i=1}^n \widehat{\gamma}_i\norm{y_i^{k} - y_i^\star}^2 \Big]
		\in \mathcal{O}((1 - \rho)^k).
	\end{align*}
	Taking the limit as $\delta \to \infty$ results in the stated interval.
\end{proof}

\begin{proof}[\textbf{Proof of \Cref{thm:pvrsg-conv-coherent}}]%
	The proof is analogous to the proof of \Cref{thm:pvrsg-conv} but with \Cref{lem:dual-contract-uni} instead of \Cref{lem:dual-contract}, yielding
	\begin{align*}
		\rho_P &= \mu\lambda(2 - \nu \lambda) \\
		\rho_D &= \min_i
		\begin{cases}
			{\textstyle \eta - (1 - \frac{np_i}{L_i}\mu )\tfrac{1}{\gamma_i}} & \mathrm{if}\text{ } \gamma_i > 0\\
			1 & \mathrm{if}\text{ } \gamma_i = 0
		\end{cases} \\
		\nu &= \max_i  {\textstyle(1+\delta^{-1})\frac{L_i\eta\gamma_i}{n p_i}  +  (1+\delta)\frac{L_i}{n p_i} - \delta\mu}.
	\end{align*}
	where $\delta > 0$, $\gamma_i \geq 0$ and $\gamma_i=0$ implies $\frac{L_i}{np_i} \leq \mu$ for all $i\in\{1,\dots,n\}$.
	Let $\gamma_i = \frac{1}{\eta - \rho_D}\max(0,1 - \frac{np_i\mu}{L_i})$ and $\delta = \sqrt{\frac{\eta}{\eta - \rho_D}}$.
	Both choices are valid if $\rho_D < \eta$ since then $\delta > 0$, $\gamma_i \geq 0$ and $\gamma_i = 0$ only if $\frac{L_i}{np_i} \leq \mu$.

	Assuming $\max_i \frac{L_i}{np_i} > \mu$ yields
	\begin{align*}
		\nu
		&= \max_i  {\textstyle(1+\delta^{-1})\frac{L_i}{n p_i}\frac{\eta}{\eta - \rho_D}\max(0,1 - \frac{np_i\mu}{L_i})  +  (1+\delta)\frac{L_i}{n p_i} - \delta\mu} \\
		&= \max_i  {\textstyle(1+\delta^{-1})\frac{\eta}{\eta - \rho_D}\max(0,\frac{L_i}{n p_i} - \mu)  +  (1+\delta)\frac{L_i}{n p_i} - \delta\mu} \\
		&= {\textstyle(1+\delta^{-1})\frac{\eta}{\eta - \rho_D}(\max_i \frac{L_i}{n p_i} - \mu)  +  (1+\delta)(\max_i\frac{L_i}{n p_i}) - \delta\mu} \\
		&= \mu +  \Big( \max_i {\textstyle \frac{L_i}{n p_i} - \mu\Big)\left(1+\sqrt{ \frac{\eta}{\eta - \rho_D}  }\right)^2}.
	\end{align*}
	Restricting the problem to $\rho = \rho_D = \rho_P$ and only considering the convergent rates $\rho \in (0,1]$ yield the problem in the theorem.
	The first statement of the theorem comes from \Cref{prop:main-conv} with $\gamma_i$ and $\delta$ inserted in the expression for $\widehat{\gamma}_i$ from \Cref{def:lyap-terms}.

	When $\max_i \frac{L_i}{np_i} = \mu$, \Cref{lem:sampled-lipschitz} gives $\frac{L_i}{np_i} = \mu$ for all $i\in\{1,\dots,n\}$, meaning $\gamma_i = 0$ is a valid choice for all $i\in\{1,\dots,n\}$.
	With this choice, $\nu = \mu$ regardless of $\delta$, and $\rho_D$ is no longer limited by $\eta$ with $\rho_D = 1$.
	The statement of the theorem then follows.
\end{proof}

\section{Proof of Corollaries}\label{sec:proof-cor}

\begin{proof}[\textbf{Proof of \Cref{cor:saga}}]%
	The expected update frequency is $\eta_i = p_i$.
	Assuming $\max_i \frac{L_i}{np_i} > \mu$ and using \Cref{thm:pvrsg-conv} the convergence rate for SAGA is given by the $\rho \in (0, p_{\min})$ that satisfies
	\begin{equation}\label{eq:saga-cor-base-prob}
		\begin{aligned}
			\rho &= \mu\lambda(2 - \nu\lambda )\\
			\nu &= \mu\min_{\delta > 0} \max_i {\textstyle(1+\delta^{-1})\frac{L_i}{n p_i \mu}\frac{p_i}{p_i - \rho}  +  (1+\delta)\frac{L_i}{n p_i \mu} - \delta}.
		\end{aligned}
	\end{equation}
	If we write $\nu$ as a function of $\rho$, this can equivalently be written as finding $\rho \in (0,p_{\min})$ such that $\rho + \lambda^2\mu \nu(\rho) = 2\mu\lambda$.
	Since $\nu(\rho)$ is continuous and $\nu(\rho) \to \infty$ as $\rho \to p_{\min}$ from below, if we find a $\tilde{\rho} \in (0,p_{\min})$ such that $\tilde{\rho} + \lambda2\mu\nu(\tilde{\rho}) \leq 2\mu\lambda$, it exists $\rho \in [\tilde{\rho},p_{\min})$ such that \cref{eq:saga-cor-base-prob} hold.
	Hence, if we replace $\nu$ in \cref{eq:saga-cor-base-prob} with an upper bound, we can find a lower bound on the contraction $\rho$.

	Let $\kappa_{\max} = \max_i \frac{L_i}{np_i\mu}$ and $p_{\min} = \min_i p_i$ and upper bound $\nu$ as
	\begin{align*}
		\nu
		&\leq \mu\min_{\delta > 0}  {\textstyle(1+\delta^{-1})\kappa_{\max}\frac{p_{\min}}{p_{\min} - \rho_D}  +  (1+\delta)\kappa_{\max} - \delta} \\
		&= \mu + \mu {\textstyle \Big[ \sqrt{\kappa_{\max} \frac{p_{\min}}{p_{\min} - \rho}} + \sqrt{\kappa_{\max} - 1} \Big]^2} \\
		&= {\textstyle \mu\kappa_{\max} \frac{2p_{\min} - \rho}{p_{\min} - \rho} + 2\mu\sqrt{\kappa_{\max}^2- \kappa_{\max}}\sqrt{\frac{p_{\min}}{p_{\min} - \rho}}} \\
		&\leq {\textstyle\mu\kappa_{\max} \frac{2p_{\min} - \rho}{p_{\min} - \rho} + \mu\sqrt{\kappa_{\max}^2 - \kappa_{\max}}\frac{2p_{\min} - \rho}{p_{\min} - \rho}} \\
		&= {\textstyle \mu\kappa_{\max}(1 + \sqrt{1 - \kappa_{\max}^{-1}})\frac{2p_{\min} - \rho}{p_{\min} - \rho}}.
	\end{align*}
	The last inequality is given by $2\sqrt{a} \leq 1 + a$ for all  $a \geq 0$.
	It can be verified that this upper bound also is valid when $\max_i \frac{L_i}{np_i} = \mu$.
	Replace $\nu$ in \Cref{eq:saga-cor-base-prob} with this upper bound gives a set of equations that define a lower bound on the contraction $\rho$.

	Inserting the two samplings and solving for $\lambda$ when the lower bound on $\rho$ is zero gives the $\lambda_{\max}$.
	Maximizing the lower bound on $\rho$ w.r.t. $\lambda$ yield the optimal $\lambda^\star$ and $\rho^\star$.
	For both uniform and Lipschitz sampling, the upper bound on $\nu$ is tight for $\rho = 0$ so it can be used to accurately determine maximal step-size according to \Cref{thm:pvrsg-conv}.
\end{proof}

\begin{proof}[\textbf{Proof of \Cref{cor:saga-sampling}}]%
	The proof is similar to the proof of \Cref{cor:saga} but instead the following upper bound is used:
	\begin{gather*}
		\nu
		\leq \mu\max_i {\textstyle 2\frac{L_i}{n p_i \mu}\frac{p_i}{p_i - \rho}  +  2\frac{L_i}{n p_i \mu} - \delta}
		\leq \mu\max_i {\textstyle 2\frac{L_i}{n p_i \mu}\frac{p_i}{p_i - \rho}  +  2\frac{L_i}{n p_i \mu}} \\
		= 2\mu \max_i {\textstyle \frac{L_i}{n\mu} \big[ \frac{1}{p_i} + \frac{1}{p_i - \rho} \big]}.
	\end{gather*}
	Replacing $\nu$ in \Cref{thm:pvrsg-conv} with this upper bound and inserting the presented $p_i$, $\lambda$ and $\rho$ verifies the first claim.

	The rate from \Cref{thm:pvrsg-conv} is of the form $\Exp \norm{x^k - x^\star}^2 \in \mathcal{O}( (1 - \lambda^\star\mu)^k )$.
	The iteration complexity to achieve an $\epsilon$-accurate solution in expectation is then $k \in \mathcal{O}(\frac{1}{\lambda^\star\mu}\log \frac{1}{\epsilon})$.
	One gradient evaluation is done per iteration so $\mathcal{O}(\frac{1}{\lambda^\star\mu}\log \frac{1}{\epsilon})$ is also the computational complexity.
	Inserting $ \lambda^\star$ gives the result.
\end{proof}

\begin{proof}[\textbf{Proof of \Cref{cor:lsvrg}}]%
	The proof is analogous to \Cref{cor:saga} but \Cref{thm:pvrsg-conv-coherent} is used instead of \Cref{thm:pvrsg-conv} and $\nu$ is upper bounded by
	\begin{gather*}
		\nu
		\leq \mu + 2\mu( \kappa_{\max} - 1){\textstyle \left[\frac{\eta}{\eta - \rho} + 1 \right]}
		\leq \frac{\mu}{2}{\textstyle \left[\frac{\eta}{\eta - \rho} + 1 \right]} + 2\mu( \kappa_{\max} - 1){\textstyle \left[\frac{\eta}{\eta - \rho} + 1 \right]} \\
		\leq {\textstyle \mu (2\kappa_{\max} - \frac{3}{2}) \frac{2\eta - \rho}{\eta - \rho}}
	\end{gather*}
	where $\kappa_{\max} = \max_i \frac{L_i}{np_i\mu}$.
\end{proof}

\begin{proof}[\textbf{Proof of \Cref{cor:lsvrg-complexity-full}}]%
	From \Cref{cor:lsvrg} we get the iteration complexity $k \in \mathcal{O}(\frac{1}{\lambda^\star\mu}\log \frac{1}{\epsilon})$.
	One gradient evaluation is needed for the primal update and $n\eta$ evaluations are needed in expectation for the dual update, this gives the computational complexity $\mathcal{O}((1+n\eta)\frac{1}{\lambda^\star \mu}\log \frac{1}{\epsilon})$.
	Inserting $\lambda^\star$ from \Cref{cor:lsvrg} and using $\tfrac{1}{2}(a + b + \sqrt{a^2 + b^2}) \leq a + b$ gives the result.
\end{proof}

\end{appendices}

\bibliography{references}

\end{document}